# Estimation of bivariate excess probabilities for elliptical models

BELKACEM ABDOUS[1], ANNE-LAURE FOUGÈRES[2,*],
KILANI GHOUDI[3] and PHILIPPE SOULIER[2,**]

[1]*Département de médecine sociale et préventive, Université Laval, Québec, Canada G1K 7P4.
E-mail: belkacem.abdous@msp.ulaval.ca*

[2]*Université Paris Ouest-Nanterre, Equipe MODAL'X, 92000 Nanterre, France.
E-mail: [*]Anne-Laure.Fougeres@u-paris10.fr; [**]Philippe.Soulier@u-paris10.fr*

[3]*Department of Statistics, United Arab Emirates University, P.O Box 17555 Al-Ain, United Arab Emirates. E-mail: kghoudi@uaeu.ac.ae*

Let $(X,Y)$ be a random vector whose conditional excess probability $\theta(x,y) := P(Y \leq y \mid X > x)$ is of interest. Estimating this kind of probability is a delicate problem as soon as $x$ tends to be large, since the conditioning event becomes an extreme set. Assume that $(X,Y)$ is elliptically distributed, with a rapidly varying radial component. In this paper, three statistical procedures are proposed to estimate $\theta(x,y)$ for fixed $x,y$, with $x$ large. They respectively make use of an approximation result of Abdous *et al.* (cf. *Canad. J. Statist.* **33** (2005) 317–334, Theorem 1), a new second order refinement of Abdous *et al.*'s Theorem 1, and a non-approximating method. The estimation of the conditional quantile function $\theta(x, \cdot)^{\leftarrow}$ for large fixed $x$ is also addressed and these methods are compared via simulations. An illustration in the financial context is also given.

## 1. Introduction

Consider two positively dependent market returns, $X$ and $Y$. It is of practical importance to assess the possible contagion between $X$ and $Y$. Contagion formalizes the fact that for large values $x$, the probability $P(Y > y \mid X > x)$ is greater than $P(Y > y)$: see, for example, Abdous *et al.* [1] or Bradley and Taqqu [9, 10], among others. Besides, this conditional probability is also related to the tail dependence coefficient which has been widely investigated in the financial risk management context: see, for instance, Frahm *et al.* [22]. Therefore, the behavior of the conditional excess probability $\theta(x,y) := P(Y \leq y \mid X > x)$ is of practical interest, especially for large values of $x$. Estimating this kind of probability is a delicate problem as soon as $x$ tends to be large, since the conditioning event becomes an extreme set. "Large" here essentially means that $x$ is beyond the largest value of the $X$ observations so that the conditional empirical distribution function then fails to be of any use, even if the probability $\theta(x, y)$ in itself is *not* a small probability, nor close to 1. Alternative methods have to be considered.

A classical approach is to call on multivariate extreme value theory. Many refined inference procedures have been developed, making use of the structure of multivariate max-stable distributions introduced by De Haan and Resnick [15], Pickands [37] and De Haan [14]. These procedures are successful in the rather general situation where $(X, Y)$ are asymptotically dependent (for the maxima), which means heuristically that $X$ and $Y$ can be simultaneously large (see, e.g., Resnick [38] or Beirlant *et al.* [5] for more details).

Efforts have recently been made to address the problem in the opposite case of asymptotic independence. In some papers, attempts are made to provide models for joint tails (see, e.g., Ledford and Tawn [33, 34], Draisma *et al.* [18], Resnick [39], Maulik and Resnick [36]). In a parallel way, Heffernan and Tawn [28] explored modeling for multivariate tail regions which are not necessarily joint tails, and Heffernan and Resnick [27] provided a complementary mathematical framework in the bivariate case.

In these papers, the key assumption is that there exists a limit for the conditional distribution of $Y$, suitably centered and renormalized, given that $X$ tends to infinity. This assumption was first checked for bivariate spherical distributions by Eddy and Gale [19] and Berman [8], Theorem 12.4.1. Abdous *et al.* [1] obtained it for bivariate elliptical distributions, Hashorva [25] for multivariate elliptical distributions, Balkema and Embrechts [3] for generalized multivariate elliptical distributions and Hashorva *et al.* [26] for Dirichlet multivariate distributions.

Elliptical distributions form a large family of multivariate laws which have received considerable attention, especially in the financial risk context; see Artzner *et al.* [2], Embrechts *et al.* [20], Hult and Lindskog [30], among others. Assume from now on that $(X, Y)$ is elliptically distributed; Theorem 1 of Abdous *et al.* [1] exhibits the asymptotic behavior of $\theta(x, y)$ when $x \to \infty$ for such an elliptical pair. The appropriate rate $y = y(x)$ is made explicit to get a non-degenerate behavior of $\lim_{x \to \infty} \theta(x, y)$. This rate depends on the tail behavior of the radial random variable $R$ defined by the relation $R^2 = (X^2 - 2\rho XY + Y^2)/(1 - \rho^2)$, where $\rho$ is the Pearson correlation coefficient between $X$ and $Y$. The only parameters involved in $y(x)$ and in the limiting distribution are the Pearson correlation coefficient $\rho$ and the index of regular variation of $R$ (say $\alpha$) or an auxiliary function of $R$ (say $\psi$), depending on whether $R$ has a regularly or rapidly varying upper tail.

In financial applications, the regularly varying behavior of the upper tails is commonly encountered. Abdous *et al.* [1] provided a simulation study in the specific case where $R$ has a regularly varying tail. Existing estimators of $\rho$ and $\alpha$ were used therein to obtain a practical way to estimate excess probabilities. However, this assumption fails to hold in some situations, as shown by Levy and Duchin [35] who compared the fit of 11 distributions on a wide range of stock returns and investment horizons. They concluded that the logistic distribution, which has a rapidly varying tail, gives the best fit in most of the cases for weekly and monthly returns.

The aim of this paper is to focus on the case where the radial component $R$ associated with the elliptical pair $(X, Y)$ has a rapidly varying upper tail. A second order approximation result is obtained, which refines Theorem 1 of Abdous *et al.* in the case of rapid variation of $R$. As we mentioned earlier, under the rapidly varying upper tail assumption, the obtained approximation involves an auxiliary function $\psi$ which has to be estimated. To the best of our knowledge, estimation of the auxiliary function has never been considered in the literature. We propose three statistical procedures for the estimation of $\theta(x, y)$ for fixed $x, y$ with $x$ large. They respectively make use of Abdous *et al.*'s approximation result, its second order refinement and an alternative method which does not rely on an analytic approximation result.

The conditional quantiles – also called *regression quantiles* – of a response $Y$ given a covariate $X$ have received significant attention; see, for example, Koenker and Bassett [32] or Koenker [31]. More specifically, extremal conditional quantiles have proven to be useful in economics and financial applications (see, e.g., Chernozhukov [11]). In the same spirit, one might be interested in estimating the conditional quantile function $\theta(x, \cdot)^{\leftarrow}$ when $x$ is fixed and large. This problem is also addressed in detail in the present paper and the simulation study performed is presented in terms of estimation of both $\theta$ and $\theta(x, \cdot)^{\leftarrow}$.

The paper is organized as follows. The second order approximation result is presented in Section 2, as well as some remarks and examples illustrating the theorem. The statistical procedures are described in Section 3, where a semi-parametric estimator of $\psi$ is proposed. Section 4 deals with a comparative simulation study, while Section 5 provides an application to the financial context, revisiting data studied by Levy and Duchin [35]. Some concluding comments are given in Section 6. Proofs are deferred to the Appendix.

## 2. Asymptotic approximation

Consider a bivariate elliptical random vector $(X, Y)$. General background on elliptical distributions can be found in, for example, Fang *et al.* [21]. One can focus, without loss of theoretical generality, on the standard case where $\mathbb{E}X = \mathbb{E}Y = 0$ and $\operatorname{Var} X = \operatorname{Var} Y = 1$. A convenient representation is then the following (see, e.g., Hult and Lindskog [30]): $(X, Y)$ has a standard elliptical distribution with radial positive distribution function $H$ and Pearson correlation coefficient $\rho$ if it can be expressed as

$$(X, Y) = R(\cos U, \rho \cos U + \sqrt{1 - \rho^2} \sin U),$$

where $R$ and $U$ are independent, $R$ has distribution function $H$ with $\mathbb{E}R^2 = 2$ and $U$ is uniformly distributed on $[-\pi/2, 3\pi/2]$. Hereafter, to avoid trivialities, we assume $|\rho| < 1$.

Let $\Phi$ denote the normal distribution function and $\varphi$ its density, that is,

$$\varphi(t) = \frac{\mathrm{e}^{-t^2/2}}{\sqrt{2\pi}}, \qquad \Phi(x) = \int_{-\infty}^{x} \varphi(t) \, \mathrm{d}t.$$

This paper deals with elliptical distributions with rapidly varying marginal upper tails or, equivalently, with a rapidly varying radial component. More precisely, the radial component $R$ associated with $(X, Y)$ is assumed to be such that there exists an *auxiliary function* $\psi$ for which one gets, for any positive $t$,

$$\lim_{x \to \infty} \frac{P\{R > x + t\psi(x)\}}{P(R > x)} = e^{-t}. \tag{1}$$

Such a condition implies that $R$ belongs to the max-domain of attraction of the Gumbel distribution; see Resnick [38], page 26 for more details. De Haan [13] introduced this class of distributions as the $\Gamma$-*varying class*. The function $\psi$ is positive, absolutely continuous and satisfies $\lim_{t \to \infty} \psi'(t) = 0$, $\lim_{t \to \infty} \psi(t)/t = 0$ and $\lim_{t \to \infty} \psi\{t + x\psi(t)\}/\psi(t) = 1$ for each positive $x$. It is only unique up to asymptotic equivalence.

Let us recall Abdous *et al.* [1]'s result in this rapidly varying context (see Theorem 1(ii)).

**Theorem 1.** *Let $(X, Y)$ be a bivariate standardized elliptical random variable with Pearson correlation coefficient $\rho$ and radial component satisfying (1). Then, for each $z \in \mathbb{R}$, one has*

$$\lim_{x \to \infty} P(Y \le \rho x + z\sqrt{1-\rho^2}\sqrt{x\psi(x)} \mid X > x) = \Phi(z).$$

In the following subsection, a rate of convergence is provided for the approximation result stated in Theorem 1, as well as a second order correction. Note that if $(X, Y)$ has an elliptic distribution with correlation coefficient $\rho$, then the couple $(X, -Y)$ has an elliptic distribution with correlation coefficient $-\rho$. Therefore, one can focus on nonnegative $\rho$. From now on, assume that $\rho \ge 0$. As a consequence, one can assume that both $x > 0$ and $y > 0$.

## 2.1. Main result

For each distribution function $H$, denote by $\bar{H}$ the survival function $\bar{H} = 1 - H$. The following assumption will be sufficient to obtain the main result. It is a strengthening of (1).

*Assumption 1.* *Let $H$ be a rapidly varying distribution function such that*

$$\left| \frac{\bar{H}\{x + t\psi(x)\}}{\bar{H}(x)} - e^{-t} \right| \le \chi(x)\Theta(t) \tag{2}$$

*for all $t \ge 0$ and $x$ large enough, where $\lim_{x \to \infty} \chi(x) = 0$, $\psi$ satisfies*

$$\lim_{x \to \infty} \frac{\psi(x)}{x} = 0 \tag{3}$$

*and $\Theta$ is locally bounded and integrable over $[0, \infty)$.*

The following result is a second order approximation for conditional excess probabilities in the elliptical case with rapidly varying radial component.

**Theorem 2.** *Let $(X, Y)$ be a bivariate elliptical vector with Pearson correlation coefficient $\rho \in [0, 1)$ and radial distribution $H$ that satisfies Hypothesis 1. Then, for all $x > 0$ and $z \in \mathbb{R}$,*

$$\mathbb{P}(Y \leq \rho x + z\sqrt{1-\rho^2}\sqrt{x\psi(x)} \mid X > x)$$
$$= \Phi(z) - \sqrt{\frac{\psi(x)}{x}} \frac{\rho \varphi(z)}{\sqrt{1-\rho^2}} + O\left(\chi(x) + \frac{\psi(x)}{x}\right), \quad (4)$$

$$\mathbb{P}(Y \leq \rho x + z\sqrt{1-\rho^2}\sqrt{x\psi(x)} + \rho\psi(x) \mid X > x)$$
$$= \Phi(z) + O\left(\chi(x) + \frac{\psi(x)}{x}\right). \quad (5)$$

*All the terms $O(\ )$ are locally uniform with respect to $z$.*

**Remark 1.** This result provides a rate of convergence in the approximation result of Abdous *et al.* [1] Theorem 1, and a second order correction. This correction is useful only if $\chi(x) = o(\sqrt{\psi(x)/x})$.

**Remark 2.** Theorem 2 and the formula

$$\mathbb{P}(X \leq x'; Y \leq y \mid X > x)$$
$$= \mathbb{P}(Y \leq y \mid X > x) - \mathbb{P}(X > x' \mid X > x)\mathbb{P}(Y \leq y \mid X > x')$$

yield (with some extra calculations) the asymptotic joint distribution (and the rate of convergence) of $(X, Y)$ given that $X > x$ when $x$ is large: for all $x > 0$ and $z \in \mathbb{R}$,

$$\mathbb{P}(X \leq x + t\psi(x); Y \leq \rho x + z\sqrt{1-\rho^2}\sqrt{x\psi(x)} \mid X > x)$$
$$= (1 - \mathrm{e}^{-t})\Phi(z) + O\left(\chi(x) + \sqrt{\frac{\psi(x)}{x}}\right), \quad (6)$$

where all the terms $O(\ )$ are locally uniform with respect to $z$.

**Remark 3.** Hashorva [25] obtained that the conditional limit distribution of $Y$, given that $X = x$ for multivariate elliptical vectors with rapidly varying radial component, is also the Gaussian distribution. The equality of these two asymptotic distributions is not true in general. In the elliptical context, this is only true if the radial variable is rapidly varying. The conditional distribution of $Y$ given that $X = x$ is of course related to the joint distribution of $(X, Y)$ given that $X > x$ via the formula

$$\mathbb{P}(X \leq x'; Y \leq y \mid X > x) = \int_x^{x'} \mathbb{P}(Y \leq y \mid X = u) P_X(\mathrm{d}u),$$

where $P_X$ is the distribution of $X$. However, the limiting behavior of the integrand is not sufficient to obtain the limit of the integral, so (6) is not a straightforward consequence of Hashorva's result.

## 2.2. Remarks and examples

Hypothesis 1 gives a rate of convergence in the conditional excess probability approximation. To the best of our knowledge, the literature deals more classically with second order conditions providing limits (see, e.g., Beirlant *et al.* [5], Section 3.3) or with pointwise or uniform rates of convergence (cf. De Haan and Stadtmüller [16] or Beirlant *et al.* [6]). The need here is to have a non-uniform bound that can be used for dominated convergence arguments. However, in the examples given below, the non-uniform rates $\chi(x)$ that we exhibit are the same as the optimal uniform rates provided by De Haan and Stadtmüller [16].

One can, however, check that Hypothesis 1 holds for the usual rapidly varying functions, in particular, for most of the so-called *Von Mises distribution functions*, which satisfy (see, e.g., Resnick [38], page 40)

$$\bar{H}(x) = d \exp\left\{-\int_{x_0}^x \frac{\mathrm{d}s}{\psi_H(s)}\right\} \tag{7}$$

for $x$ greater than some $x_0 \geq 0$, where $d > 0$ and $\psi_H$ is positive, absolutely continuous and $\lim_{x\to\infty} \psi_H'(x) = 0$. Note that under this assumption, $\psi_H = \bar{H}/H'$ and $\psi_H$ is an auxiliary function in the sense of (1). In the sequel, auxiliary functions $\psi_H$ satisfying (7) will be called *Von Mises auxiliary functions*.

The following lemma provides sufficient conditions for Hypothesis 1, which could be weakened at the price of additional technicalities.

**Lemma 1.** *Let $H$ be a Von Mises distribution function and assume that*

(i) *$\psi_H$ is ultimately monotone, differentiable and $|\psi_H'|$ is ultimately decreasing;*

(ii) *if $\psi_H$ is decreasing, then either $\lim_{x\to\infty} \psi_H(x) > 0$ or there exist positive constants $c_1$ and $c_2$ such that for all $x \geq 0$ and $u \geq 0$,*

$$\frac{\psi_H(x)}{\psi_H(x+u)} \leq c_1 \mathrm{e}^{c_2 u}. \tag{8}$$

*Then Hypothesis 1 holds the for $\psi_H$:*

$$\left|\frac{\bar{H}\{x + t\psi_H(x)\}}{\bar{H}(x)} - \mathrm{e}^{-t}\right| \leq \chi(x)\Theta(t),$$

*with $\chi(x) = O(|\psi_H'(x)|)$ and $\Theta(t) = O((1+t)^{-\kappa})$ for an arbitrary $\kappa > 0$.*

**Remark 4.** Assumption (8) holds if $\psi_H$ is regularly varying with index $\gamma < 0$. It also holds if $\psi_H(x) = \mathrm{e}^{-cx}$ for some $c > 0$.

***Remark 5.*** If $\psi_H$ is regularly varying with index $\gamma < 1$ and $\psi'_H$ ultimately decreasing, then $\psi'_H(x) = o(\sqrt{\psi_H(x)/x})$ and hence the second order correction is useful.

***Remark 6.*** This bound corresponds to a worst-case scenario. In many particular cases, a much faster rate of convergence can be obtained. For instance, if $\bar{H}(x) = e^{-t}$, then (2) holds with $\psi_H \equiv 1$, whence $\chi(x) = 1/x$, but for any positive $x$ and $t$, $\bar{H}\{x + t\psi_H(x)\}/\bar{H}(x) = e^{-t}$. The rate of convergence is infinite here. If $\bar{H}$ is the Gumbel distribution, then $\psi_H \equiv 1$ and the rate of convergence in (2) is exponential: $|\bar{H}(x+t)/\bar{H}(x) - e^{-t}| \leq 2e^{-x}e^{-t}$.

Some examples of continuous distributions satisfying Hypothesis 2 are given for the purposes of illustration in Section 4, Table 1. The following example illustrates a discrete situation.

***Example 1 (Discrete distribution in the domain of attraction of the Gumbel law).*** Let $\psi$ be a concave increasing function such that $\lim_{x \to \infty} \psi(x) = +\infty$ and $\lim_{x \to \infty} \psi'(x) = 0$, and define

$$\bar{H}_\#(x) = \exp\left\{-\int_{x_0}^{x} \frac{ds}{\psi(s)}\right\}, \qquad \bar{H}(x) = \bar{H}_\#([x]),$$

where $[x]$ is the integer part of $x$. $H$ is then a discrete distribution function and $\bar{H}$ belongs to the domain of attraction of the Gumbel distribution, but does not satisfy condition (7). Nevertheless, following the lines of the proof of Lemma 1, one can check that condition (2) holds with

$$\chi(x) = O(\psi'([x]) + 1/\psi(x)), \qquad \Theta(t) = O((1+t)^{-\kappa})$$

for any arbitrary $\kappa > 0$.

One can deduce from Lemma 1 that Theorem 2 holds for $\psi = \psi_H$. The following lemma concerns what happens if one uses an asymptotically equivalent auxiliary function $\psi$ instead of $\psi_H$ in Theorem 2.

**Lemma 2.** *Under the assumption of Lemma 1, let $\psi$ be equivalent to $\psi_H$ at infinity and define $\xi(x) = |\psi(x) - \psi_H(x)|/\psi_H(x)$. Then*

$$\left|\frac{\bar{H}(x + t\psi(x))}{\bar{H}(x)} - e^{-t}\right| \leq O\{|\psi'_H(x)| + \xi(x)\}\Theta(t). \tag{9}$$

***Remark 7.*** A consequence of Lemma 2 is that if an auxiliary function $\psi$ is used instead of $\psi_H$ in Theorem 2, then the second order correction is relevant, provided that $\xi(x) = o(\sqrt{\psi_H(x)/x})$. This is the case in the examples given on line 1 and line 3 of Table 1 ('Normal' and 'Logis') if one takes $\psi(x) = 1/x$ or $\psi(x) = 1/(2x)$, respectively, and in the example given on line 5 of Table 1 ('Lognor') when using $\psi(x) = x/\log(x)$.

## 3. Statistical procedure

For *given* large positive $x$ and $y$, consider the problems of estimating $\theta(x,y) = \mathbb{P}(Y \leq y \mid X > x)$ and the conditional quantile function $\theta(x, \cdot)^{\leftarrow}$. Note, in passing, that in a practical situation, $x$ is neither a threshold nor a parameter of the statistical procedure, but a value imposed by the practical problem (e.g., a high quantile of the marginal distribution of $X$). The empirical distribution function is useless since there might be no observations in the considered range. We suggest estimating these quantities by means of Theorem 2.

Assume that a sample $(X_1, Y_1), \ldots, (X_n, Y_n)$ is available, drawn from an elliptical distribution with radial component satisfying Hypothesis 2. Note that in this section, we do not assume that the distribution is standardized, so a preliminary standardization is required.

### 3.1. Definition of the estimators

The estimation of $\theta(x,y)$ and the conditional quantile function requires estimates of $\mu_X$, $\mu_Y$, $\sigma_X$, $\sigma_Y$, $\rho$ and $\psi$. Let $\hat{\mu}_X$, $\hat{\mu}_Y$, $\hat{\sigma}_X$, $\hat{\sigma}_Y$, $\hat{\rho}_n$ and $\hat{\psi}_n$ denote such estimates. For fixed $x, y > 0$, define

$$\hat{\theta}_{n,1}(x,y) = \Phi\left(\frac{\hat{y} - \hat{\rho}_n \hat{x}}{\sqrt{1 - \hat{\rho}_n^2}\sqrt{\hat{x}\hat{\psi}_n(\hat{x})}}\right) \tag{10}$$

and

$$\hat{\theta}_{n,2}(x,y) = \Phi\left(\frac{\hat{y} - \hat{\rho}_n \hat{x} - \hat{\rho}_n \hat{\psi}_n(\hat{x})}{\sqrt{1 - \hat{\rho}_n^2}\sqrt{\hat{x}\hat{\psi}_n(\hat{x})}}\right), \tag{11}$$

where $\hat{x} = (x - \hat{\mu}_X)/\hat{\sigma}_X$ and $\hat{y} = (y - \hat{\mu}_Y)/\hat{\sigma}_Y$.

In order to estimate the conditional quantile function $\theta(x, \cdot)^{\leftarrow}$, define, for fixed $\theta \in (0,1)$,

$$\hat{y}_{n,1} = \hat{\mu}_Y + \hat{\sigma}_Y\{\hat{\rho}_n \hat{x} + \sqrt{1 - \hat{\rho}_n^2}\sqrt{\hat{x}\hat{\psi}_n(\hat{x})}\Phi^{-1}(\theta)\}, \tag{12}$$

$$\hat{y}_{n,2} = \hat{\mu}_Y + \hat{\sigma}_Y\{\hat{\rho}_n \hat{x} + \hat{\rho}_n \hat{\psi}_n(\hat{x}) + \sqrt{1 - \hat{\rho}_n^2}\sqrt{\hat{x}\hat{\psi}_n(\hat{x})}\Phi^{-1}(\theta)\}. \tag{13}$$

Estimating $\mu_X$, $\mu_Y$, $\sigma_X$, $\sigma_Y$ and $\rho$ is a classical topic and the empirical version of each quantity can easily be used. Under the assumption of elliptical distributions, however, one can observe better stability when the Pearson correlation coefficient is estimated by $\hat{\rho}_n = \sin(\pi \hat{\tau}_n / 2)$, where $\hat{\tau}_n$ is the empirical Kendall's tau (see, e.g., Hult and Lindskog [30] for more details). This estimator is $\sqrt{n}$-consistent and asymptotically normal.

Consider now the problem of the estimation of $\psi$. Since auxiliary functions are defined up to asymptotic equivalence, a particular representantive must be a priori chosen in order to define the estimator. We assume that an admissible auxiliary function is

$$\psi(x) = \frac{1}{c\beta}x^{1-\beta} \tag{14}$$

for some constants $c > 0$ and $\beta > 0$. Under this assumption, estimation of $\psi$ reduces to estimating $c$ and $\beta$.

An extensive body of literature exists on estimators of $\beta$ (see, e.g., Beirlant *et al.* [4], Gardes and Girard [23], Dierckx *et al.* [17], among others). The method chosen here is the one proposed in Beirlant *et al.* [7]. Let $k_n$ be a user-chosen threshold and $R_{j,n}$, $1 \leq j \leq n$, be the order statistics of the sample $R_1, \ldots, R_n$. The estimator of $\beta$ is obtained as the slope of the Weibull quantile plot at the point $(\log\log(n/k_n), \log(R_{n-k_n,n}))$:

$$\hat{\beta}_n = \frac{k_n^{-1}\sum_{i=1}^{k_n} \log\log(n/i) - \log\log(n/k_n)}{k_n^{-1}\sum_{i=1}^{k_n} \log(R_{n-i+1,n}) - \log(R_{n-k_n,n})}. \tag{15}$$

An estimator of $c$ is then naturally given by

$$\hat{c}_n = \frac{1}{k_n}\sum_{i=1}^{k_n} \frac{\log(n/i)}{R_{n-i+1,n}^{\hat{\beta}_n}}. \tag{16}$$

Actually, in our context, the radial component is not observed. We estimate the $R_i$'s by

$$\hat{R}_i^2 = \hat{X}_i^2 + (\hat{Y}_i - \hat{\rho}_n\hat{X}_i)^2/(1 - \hat{\rho}_n^2),$$

where $\hat{X}_i = (X_i - \hat{\mu}_X)/\hat{\sigma}_X$ and $\hat{Y}_i = (Y_i - \hat{\mu}_Y)/\hat{\sigma}_Y$, and we plug these values into (15) and (16). We then define $\hat{\psi}_n(x) = x^{1-\hat{\beta}_n}/(\hat{c}_n\hat{\beta}_n)$.

## 3.2. Discussion of the estimation error versus approximation error

Let

$$\hat{z}_{n,1} = \frac{y - \hat{\rho}_n x}{\sqrt{1 - \hat{\rho}_n^2}\sqrt{x\hat{\psi}_n(x)}}, \qquad \hat{z}_{n,2} = \frac{y - \hat{\rho}_n x - \hat{\rho}_n\hat{\psi}_n(x)}{\sqrt{1 - \hat{\rho}_n^2}\sqrt{x\hat{\psi}_n(x)}},$$

$$z_1 = \frac{y - \rho x}{\sqrt{(1-\rho^2)}\sqrt{x^{2-\beta}/(c\beta)}}, \qquad z_2 = \frac{y - \rho x - \rho x^{1-\beta}/(c\beta)}{\sqrt{(1-\rho^2)}\sqrt{x^{2-\beta}/(c\beta)}}.$$

Then, for $i = 1, 2$,

$$\hat{\theta}_{n,i}(x,y) - \theta(x,y) = \Phi(\hat{z}_{n,i}) - \Phi(z_i) + \Phi(z_i) - \theta(x,y).$$

This shows that the estimators defined in (10) and (11) have two sources of error: the first one, $\Phi(\hat{z}_{n,i}) - \Phi(z_i)$, comes from the estimation of $\rho$, $\mu$, $\sigma$ and $\psi$, and the second one, $\Phi(z_i) - \theta(x,y)$, from the asymptotic nature of the approximations (4) and (5).

The order of magnitude of the estimation error can be measured by the rate of convergence of the estimators. In order to obtain a rate of convergence for the estimators $\hat{\beta}_n$ and $\hat{c}_n$, we assume that $H$ is a Von Mises distribution function with

$$\psi_H(x) = \frac{1}{c\beta} x^{1-\beta}\{1+t(x)\}, \tag{17}$$

where $t$ is a regularly varying function[1] with index $\eta\beta$ for some $\eta<0$. This implies that $\bar{H}(x) = \exp\{-cx^\beta[1+s(x)]\}$, where $s$ is also regularly varying with index $\eta\beta$. Under this assumption, the function $\psi$ defined in (14) is an admissible auxiliary function and Girard [24] has shown that $\hat{\beta}_n$ is $k_n^{1/2}$-consistent, for any sequence $k_n$ such that

$$k_n \to \infty, \qquad k_n^{1/2} \log^{-1}(n/k_n) \to 0, \qquad k_n^{1/2} b(\log(n/k_n)) \to 0,$$

where $b$ is regularly varying with index $\eta$; see Girard [24], Theorem 2 for details. Similarly, it can be shown that $k_n^{1/2}(\hat{c}_n - c) = O_P(1)$ under the same assumptions on the sequence $k_n$. Thus, for any $x$, $\hat{\psi}_n(x)$ is a $k_n^{1/2}$-consistent estimator of $\psi(x)$. Besides, $\hat{\rho}_n$, $\hat{\mu}_n$ and $\hat{\sigma}_n$ are $\sqrt{n}$-consistent, so the estimation error $\Phi(\hat{z}_{n,i}) - \Phi(z_i)$ is of order $k_n^{-1/2}$ in probability.

If $\psi_H$ satisfies (17), then Hypothesis 1 holds, and Theorem 2 and Lemma 2 provide a bound for the deterministic approximation error $\Phi(z_i) - \theta(x,y)$. Some easy computation shows that the second order correction is useful only if $\eta < -1/2$.

### 3.3. Discussion of an alternative method

The method described previously makes use of the asymptotic approximations of Theorem 2. It could be thought that a direct method, making use of (19) and an estimator of $\bar{H}$, would yield a better estimate of $\theta(x,y)$ since it would avoid this approximation step. Recall, however, that we specifically need to estimate the tail of the radial distribution so that a nonparametric estimator of $H$ cannot be considered. The traditional solution given by extreme value theory consists of fitting a parametric model for the tail. This will always induce an approximation error.

Nevertheless, there exists a situation in which the approximation error can be canceled: if the radial component is exactly Weibull distributed, that is, $\bar{H}(x) = \exp\{-cx^\beta\}$, then $\psi_H$ satisfies (14), so for any $x, u > 0$, (7) implies that

$$\frac{\bar{H}(xu)}{\bar{H}(x)} = \exp\left\{\int_x^{xu} \frac{\mathrm{d}s}{\psi(s)}\right\}.$$

Therefore, in this specific case, a consistent estimator of $\theta(x,y)$, say $\hat{\theta}_{n,3}(x,y)$, can be introduced which does not make use of any asymptotic expansion as in Theorem 2. Via

---

[1] A function $f$ is regularly varying at infinity with index $\alpha$ if for all $t>0$, $\lim_{x\to\infty} f(tx)/f(x) = t^\alpha$.

**Table 1.** Bivariate elliptical distributions used for the simulations, together with elliptical generator $g$, Von Mises auxiliary function $\psi_H$ (or equivalent), the function $\chi$ defined in Hypothesis 1 and the values of the parameters used (in addition to $\rho \in \{0.5, 0.9\}$)

| Bivariate law | Generator $g(u)$ | $\psi_H(x)$ | $\chi(x)$ | Parameters |
|---|---|---|---|---|
| Normal | $e^{-u/2}$ | $\dfrac{1}{x} + O\!\left(\dfrac{1}{x^3}\right)$ | $O(x^{-2})$ | |
| Kotz | $u^{\beta/2-1} e^{-u^{\beta/2}}$ | $x^{1-\beta}/\beta$ | $O(x^{-\beta})$ | $\beta \in \{1, 4\}$ |
| Logis* | $\dfrac{e^{-u}}{(1+e^{-u})^2}$ | $\dfrac{1+e^{-x^2}}{2x}$ | $O(x^{-2})$ | |
| Modified Kotz | $g_\star(u)^\ddagger$ | $\dfrac{x^{1-\beta}}{1+\beta\log x}$ | $O\!\left(\dfrac{x^{-\beta}}{\log x}\right)$ | $\beta = 3/2$ |
| Lognor** | $\dfrac{1}{u} e^{-(\log^2 u)/8}$ | $\dfrac{x}{\log x} + O\!\left(\dfrac{x}{\log^3 x}\right)$ | $O\!\left(\dfrac{1}{\log x}\right)$ | |
| Student | $\left(1+\dfrac{u}{\nu}\right)^{-(\nu+2)/2}$ | – | – | $\nu \in \{3, 20\}$ |

*"Logis" and "Lognor" refer to the elliptical distributions with logistic and lognormal generator, respectively.
**$g_\star(u) = \{(3/8)\log u + 1/2\} u^{-1/4} \exp\{-(1/2) u^{3/4} \log u\}$.

(19), we get, explicitly,

$$\hat\theta_{n,3}(x,y) = 1 - \frac{\int_{\arctan(\hat t_0)}^{\pi/2} \hat K(\hat x, \hat x, \cos(u))\,\mathrm{d}u + \int_{-\hat U_0}^{\arctan(\hat t_0)} \hat K(\hat x, \hat y, \sin(u+\hat U_0))\,\mathrm{d}u}{2\int_0^{\pi/2} \hat K(\hat x, \hat x, \cos(u))\,\mathrm{d}u}, \quad (18)$$

where

$$\hat K(x,y,v) = \exp\!\left\{\int_x^{y/v} \frac{\mathrm{d}s}{\hat\psi_n(s)}\right\},$$
$$\hat t_0 = (\hat y/\hat x - \hat\rho_n)/\sqrt{1-\hat\rho_n^2}, \qquad \hat U_0 = \arctan(\hat\rho_n/\sqrt{1-\hat\rho_n^2}),$$
$$\hat x = (x - \hat\mu_X)/\hat\sigma_X, \qquad \hat y = (y - \hat\mu_Y)/\hat\sigma_Y.$$

We included this estimator in the simulation study as a benchmark when looking at elliptical Kotz-distributed observations (see Table 1).

## 4. Simulation study

To assess the performance of the proposed estimators, we simulated several families of bivariate standard elliptical distributions. Recall that a standardized bivariate elliptical density function can be written as $f(x,y) = Cg\{(x^2 - 2\rho xy + y^2)/(1-\rho^2)\}$, where $g$ is called the *generator*, $\rho$ is the Pearson correlation coefficient and $C$ is a normalizing constant. The density of the radial component $R$ is given by $H'(r) = Krg(r^2)$, where $K$ is a normalizing constant (see, e.g., Fang *et al.* [21]).

The distributions used are presented in Table 1. The Pearson correlation coefficient will be either $\rho = 0.5$ or $\rho = 0.9$. Three of them (Normal, Kotz and Logis) are Von Mises distributions which satisfy both Hypothesis 1 and (17); in addition, the Von Mises auxiliary function of the Kotz distribution satisfies (14). The Lognor and the modified Kotz distributions satisfy Hypothesis 1 but not (14) and, finally, the bivariate Student distribution has a regularly varying radial component, so it does not satisfy any of the assumptions. These three distributions are used to explore the robustness of the proposed estimation method.

In each case, 200 samples of size 500 were simulated. Several values of $x$ were chosen, corresponding to the $(1-p)$-quantile of the marginal distribution of $X$, with $p = 10^{-3}$, $p = 10^{-4}$ and $p = 10^{-5}$. For each value of $x$, we computed (by numerical integration) the theoretical values of $y$ corresponding to the probability $\theta(x, y) = 0.05, 0.1, 0.2, \ldots, 0.8, 0.9, 0.95$. We then estimated $\theta(x, y)$ via the three proposed methods (cf. Section 3). For the estimation of the auxiliary function $\psi$, the threshold chosen corresponds to the highest 10% of the estimated $\hat{R}_i$'s. It must be noted that this choice is independent of $x$. For each fixed $x$, we also estimated the conditional quantile function $\theta(x, \cdot)^{\leftarrow}$ by both methods (12) and (13). We did not compute the estimated quantile function for Method 3 since it would involve the numerical inversion of the integrals which appear in (19). This is one advantage of Methods 1 and 2 over Method 3.

Some general features can be observed which conform to theoretical expectations. (i) First, in the given range of $x$ and $y$, there were hardly any observations, so the empirical conditional distribution function is useless. (ii) For a given probability $\theta$, the variability of the estimators slightly increases with $x$ for all underlying distributions. For a given $x$, the variability of the estimators is greater for medium values of $\theta$. (iii) The results for the Student distribution are as expected: if the degree of freedom $\nu$ is large, the estimation shows a high variability but moderate bias, while if $\nu$ is small, then the estimation is clearly inconsistent. (iv) The results for the Logis and modified Kotz distributions are similar to those for Gaussian distribution. (v) As described in Section 3.3, Method 3 is markedly better for the Kotz distribution only.

Hence, we have chosen to report only the results for the largest value of $x$ (corresponding to the $10^{-5}$-quantile of the marginal distribution of $X$) and the Normal, Kotz (with parameter $\beta = 1$ and $\beta = 4$) and Lognor distributions, for $\rho = 0.5$ and $\rho = 0.9$. Figures 1–4 illustrate the behavior of the estimators of the probability $\theta$: median, 2.5% and 97.5% quantiles of the estimation error $\hat{\theta}_{n,i}(x, y) - \theta(x, y)$ ($i = 1, 2, 3$) are shown as a function of the estimated probability. Figure 5 shows the estimated conditional quantile functions $\hat{y}_{n,i}(x, y)$ ($i = 1, 2$) and the theoretical conditional quantile function $y = \theta(x, \cdot)^{\leftarrow}$ for only three distributions and $\rho = 0.9$ because the results are much more stable as the correlation or the distribution vary. Median, 2.5% and 97.5% quantiles of the estimated conditional quantile function $\hat{y}_{n,i}(x, y)$ ($i = 1, 2$) are given as a function of the probability.

From these simulation results, one can see that the estimator of $\theta$ by Method 1 presents a systematic positive bias which, of course, induces an underestimation of the conditional quantile function. As expected, Method 2 corrects this systematic bias; the correction is better when $\rho$ is large. This is also true for the Lognormal generator, though to a lesser extent.

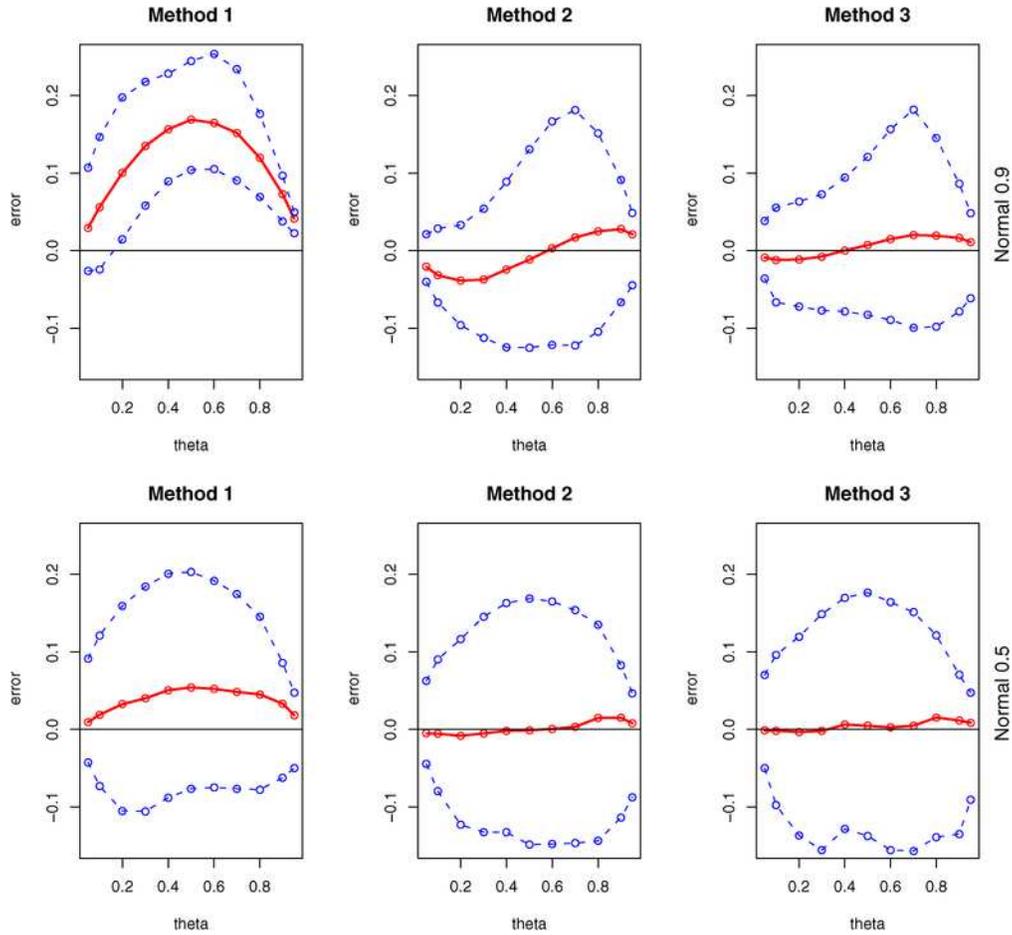

**Figure 1.** Median, 2.5% and 97.5% quantiles of the estimation error $\hat{\theta}_{n,i}(x,y) - \theta(x,y)$ ($i = 1, 2, 3$) as a function of the estimated probability. Gaussian distribution. First row: $\rho = 0.9$; second row: $\rho = 0.5$.

As already mentioned, the Lognor and modified Kotz distributions do not satisfy the assumption (17). In both cases the radial component belongs to an extended Weibull-type family, with auxiliary function $\psi$ of the form

$$\psi(x) = cx^{1-\beta}(\log x)^{-\delta}\{1 + o(1)\}$$

with $c > 0$, $\beta \geq 0$ and $\delta > 0$ if $\beta = 0$. The modified Kotz distribution corresponds to $\beta = 3/2$ and $\delta = 1$ and Lognor corresponds to $\beta = 0$ and $\delta = 1$. The simulation results are much better for the modified Kotz than for the Lognor distribution. This tends to prove that the method is not severely affected by the logarithmic factor, as long as $\beta > 0$.

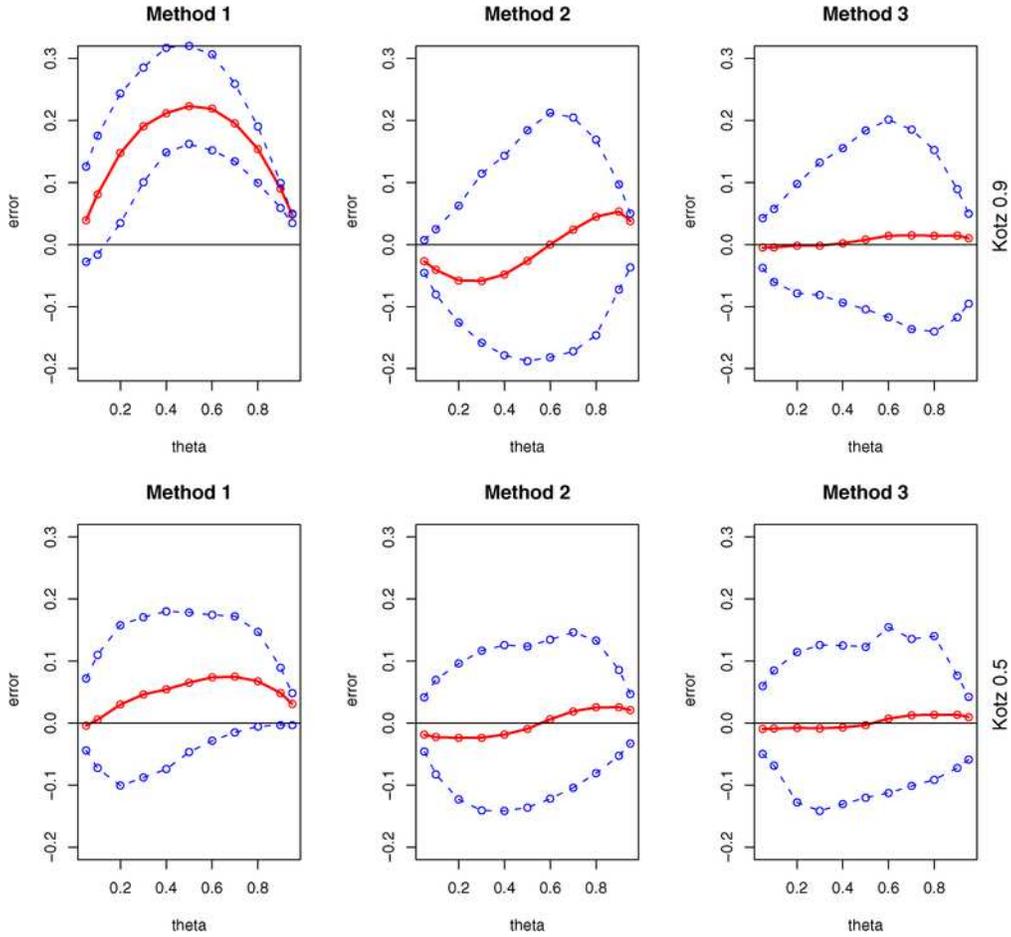

**Figure 2.** Median, 2.5% and 97.5% quantiles of the estimation error $\hat{\theta}_{n,i}(x,y) - \theta(x,y)$ ($i = 1, 2, 3$) as a function of the estimated probability. Kotz distribution, $\beta = 1$. First row: $\rho = 0.9$; second row: $\rho = 0.5$.

## 5. Financial application

In this section, the practical usefulness of our estimation procedure is illustrated in the context of financial contagion, for which an estimation of the conditional excess probability is needed. Data used by Levy and Duchin [35] are here revisited. More precisely, we consider the series of monthly returns for the 3M stock and the Dow Jones Industrial

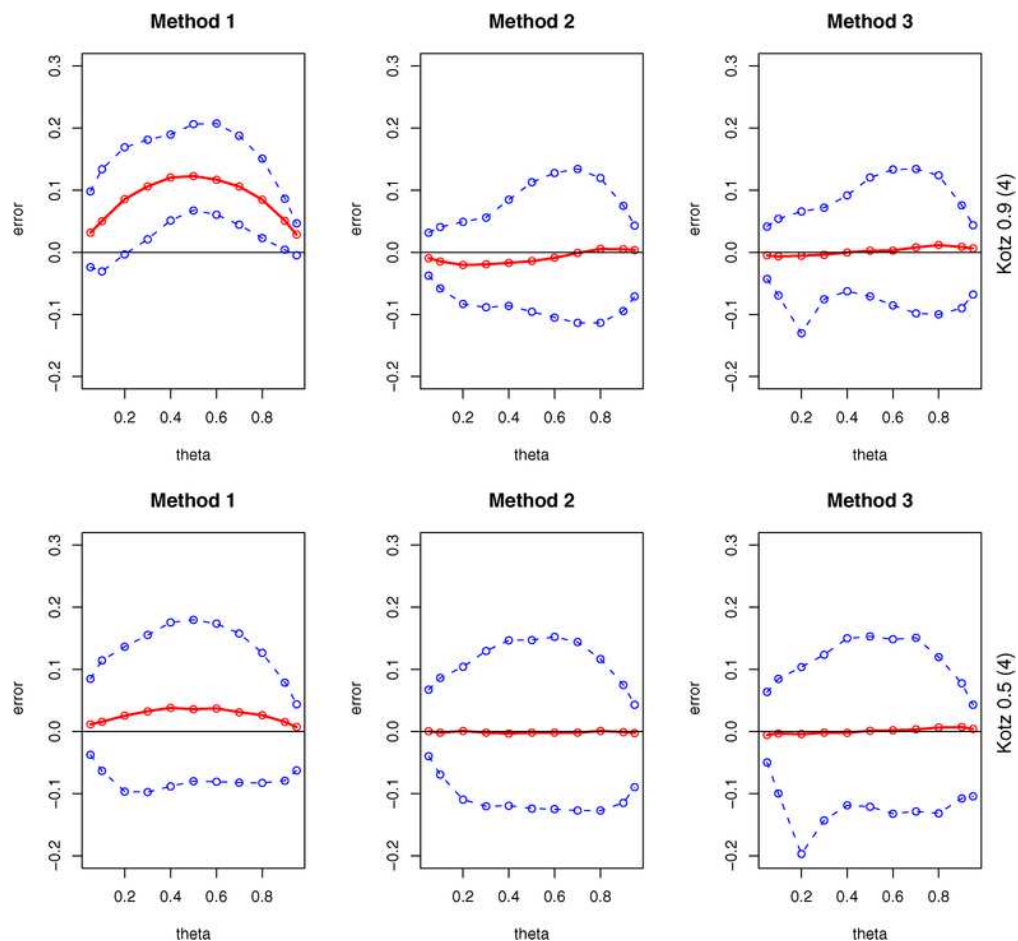

**Figure 3.** Median, 2.5% and 97.5% quantiles of the estimation error $\hat{\theta}_{n,i}(x,y) - \theta(x,y)$ ($i = 1, 2, 3$) as a function of the estimated probability. Kotz distribution, $\beta = 4$. First row: $\rho = 0.9$; second row: $\rho = 0.5$.

Average for the period ranging from January 1970 to January 2008. In the sequel, we arbitrarily investigate the conditional behavior of the 3M stock monthly returns, given some extreme values of the Dow Jones Industrial Average.

According to Levy and Duchin [35], these two series can be marginally fitted by a logistic distribution. Indeed, a Kolmogorov–Smirnov goodness-of-fit test of the logistic distribution gave us a P-value of 0.48 for the 3M returns and 0.49 for the Dow Jones Industrial Average returns. These P-values were obtained via a Monte Carlo simulation, following the procedure outlined by Stephens [40]. Moreover, the test of elliptical symmetry of Huffer and Park [29] was used to show that the data fit the bivariate elliptical model (P-value = 0.61). Finally, we checked that both marginal upper tails exhibit rapid variation; for this, we performed a generalized Pareto distribution fit to the 15%-largest values and checked via a test based on the profile likelihood 95%-confidence interval that the shape parameter could be considered equal to 0 at level 5% (see, e.g., Coles [12] for details on these classical procedures). Note that the estimated generalized Pareto tail agreed completely with the tail of the fitted logistic distribution in both cases.

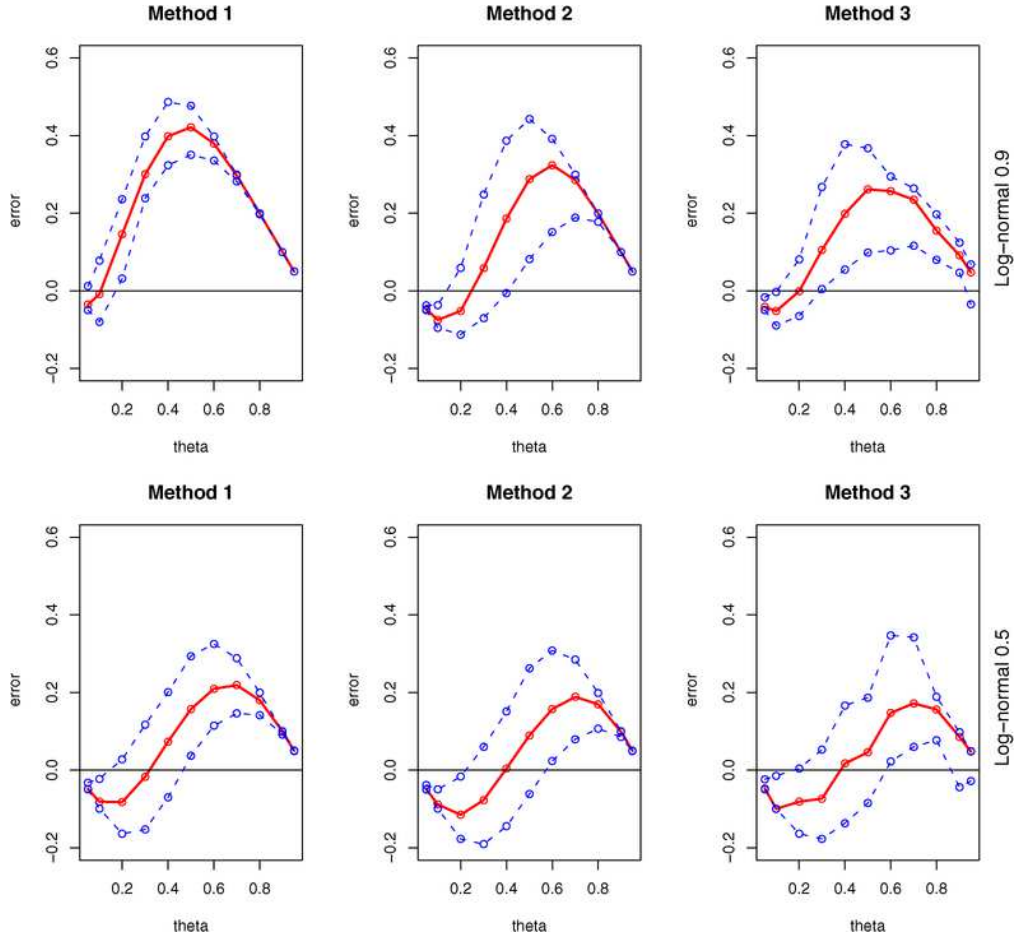

**Figure 4.** Median, 2.5% and 97.5% quantiles of the estimation error $\hat{\theta}_{n,i}(x,y) - \theta(x,y)$ ($i = 1, 2, 3$) as a function of the estimated probability. Lognor distribution. First row: $\rho = 0.9$; second row: $\rho = 0.5$.

Consequently, the estimation procedures presented in the previous sections can be applied to these data. As an illustration, in Figure 6, we depict the three estimates of $y \mapsto 1 - \theta(x,y) = \mathbb{P}(Y > y \mid X > x)$ for different values of $x$ corresponding to the 0.975, 0.99,

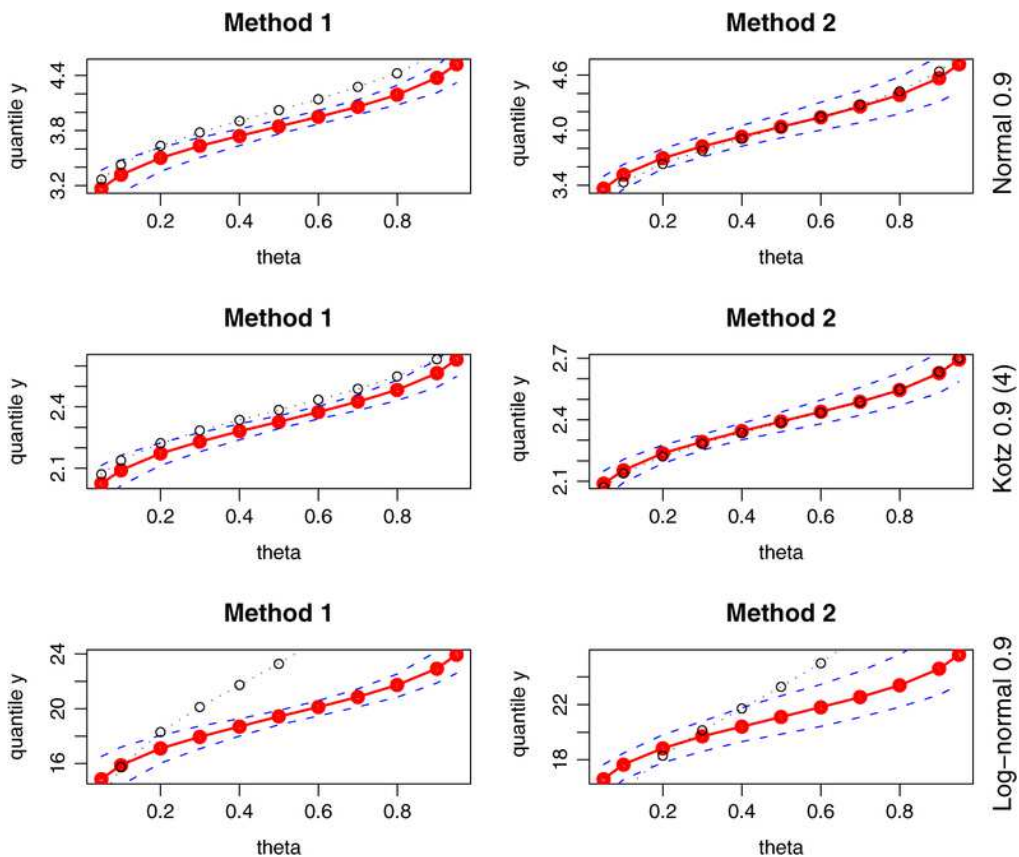

**Figure 5.** Median (solid line), 2.5% and 97.5% quantiles (dashed lines) of the estimated conditional quantile function $\hat{y}_{n,i}$ ($i = 1, 2$) defined in (12) and (13) and theoretical conditional quantile function $y$ (dotted line) as a function of the probability $\theta$. First row: Normal distribution; second row: Kotz distribution, $\beta = 4$; third row: Lognor distribution. For each of these, $\rho = 0.9$.

0.999 and 0.9999 quantiles of the fitted logistic distribution, together with the estimated marginal survival function $\mathbb{P}(Y > y)$. This last probability was estimated via the logistic distribution fitted to the $Y_i$'s. It is clearly evident from these graphics that $1 - \hat{\theta}_{n,2}$ and $1 - \hat{\theta}_{n,3}$ provide very similar estimates, whereas $1 - \hat{\theta}_{n,1}$ gives uniformly smaller values. All these estimates are uniformly greater than the marginal survival function of $Y$. This allows us to conclude that the data exhibit contagion from the Dow Jones Industrial Average to the 3M stock.

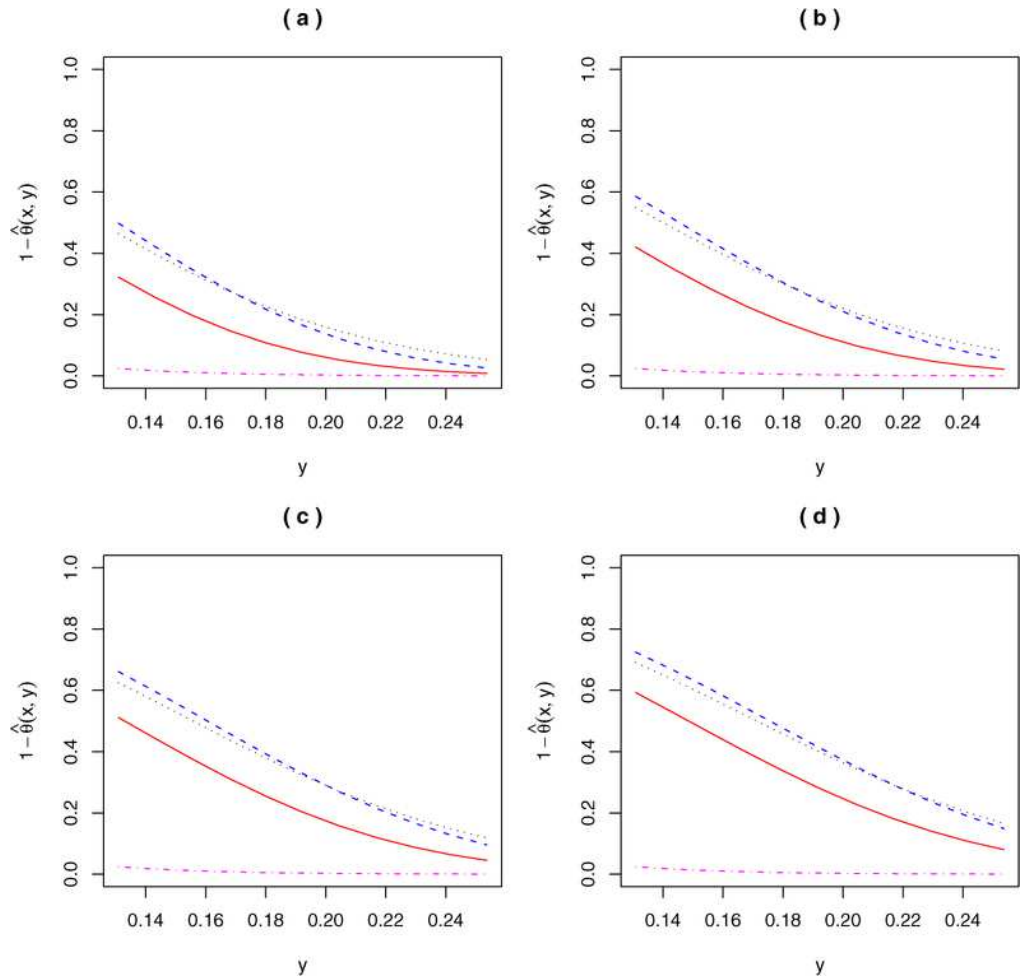

**Figure 6.** Estimates $y \mapsto 1 - \hat{\theta}_{n,i}(x,y)$ $(i=1,2,3)$ of the conditional excess distribution of the 3M stock monthly return, given that the Dow Jones Industrial Average monthly return exceeds four extreme values. The subplots (a), (b), (c) and (d) are for the values $x$ corresponding to the 0.975, 0.99, 0.999 and 0.9999 quantiles of the fitted logistic distribution, respectively. The solid line is for $\hat{\theta}_{n,1}$, the dashed line for $\hat{\theta}_{n,2}$ and the dotted line for $\hat{\theta}_{n,3}$. The estimate of the marginal survival function $P(Y > y)$ is shown as a dotted-dashed line.

## 6. Concluding remarks

In this paper, we restricted our attention to elliptically distributed random pairs $(X, Y)$ having a rapidly varying radial component. Three methods have been proposed to estimate the conditional excess probability $\theta(x, y)$ for large $x$. Under this specific assumption of elliptical distributions, Methods 2 and 3 revealed comparable results and outperformed Method 1.

Methods 1 and 2 make use of an asymptotic approximation of the conditional excess distribution function by the Gaussian distribution function. As shown by Balkema and Embrechts [3], this approximation remains valid outside the family of elliptical distributions, under geometric assumptions on the level curves of the joint density function of $(X, Y)$. This suggests that these methods may be useful outside the elliptical family. This is an ongoing research project.

# Appendix

**Proof of Theorem 2.** Define $U_0 = \arctan(\rho/\sqrt{1-\rho^2})$. For $x > 0$ and $y \in (0, x)$, we have

$$\mathbb{P}(X > x, Y > y)$$
$$= \mathbb{P}\left(R > \frac{x}{\cos U} \vee \frac{y}{\rho \cos U + \sqrt{1-\rho^2} \sin U} \,;\, -U_0 \leq U \leq \frac{\pi}{2}\right)$$

Let $t_0 = (y/x - \rho)/\sqrt{1-\rho^2}$. Then $-U_0 < \arctan(t_0)$, and $x/\cos u > y/(\rho \cos u + \sqrt{1-\rho^2} \sin u)$ if and only if $u > \arctan(t_0)$. Hence,

$$\mathbb{P}(X > x, Y > y)$$
$$= \int_{\arctan(t_0)}^{\pi/2} \bar{H}\left\{\frac{x}{\cos u}\right\} \frac{\mathrm{d}u}{2\pi} + \int_{-U_0}^{\arctan(t_0)} \bar{H}\left\{\frac{y}{\sin(u+U_0)}\right\} \frac{\mathrm{d}u}{2\pi},$$

$$\mathbb{P}(Y > y \mid X > x)$$
$$= \frac{\int_{\arctan(t_0)}^{\pi/2} \bar{H}(x/\cos(u))\,\mathrm{d}u + \int_{-U_0}^{\arctan(t_0)} \bar{H}(y/\sin(u+U_0))\,\mathrm{d}u}{2\int_0^{\pi/2} \bar{H}(x/\cos(u))\,\mathrm{d}u}. \tag{19}$$

If $t_0 \geq 0$, that is, $y - \rho x \geq 0$, the changes of variables $v = 1/\cos(u)$ and $v = 1/\sin(u + U_0)$ yield

$$\mathbb{P}(Y > y \mid X > x) = \frac{I_1 + I_2}{I_3} \tag{20}$$

with

$$I_1 = \int_{w_1}^{\infty} \frac{\bar{H}(vx)}{\bar{H}(x)} \frac{\mathrm{d}v}{v\sqrt{v^2-1}},$$
$$I_2 = \int_{w_2}^{\infty} \frac{\bar{H}(vy)}{\bar{H}(x)} \frac{\mathrm{d}v}{v\sqrt{v^2-1}},$$

$$I_3 = 2\int_1^\infty \frac{\bar{H}(vx)}{\bar{H}(x)} \frac{\mathrm{d}v}{v\sqrt{v^2-1}},$$

$$w_1 = \sqrt{1+t_0^2} = \sqrt{1+(y/x-\rho)^2/(1-\rho^2)},$$

$$w_2 = xw_1/y.$$

If $t_0 < 0$, then

$$\mathbb{P}(Y > y \mid X > x) = \frac{I_3 - I_1 + I_2}{I_3}. \tag{21}$$

Denote $w_0 = x(w_1 - 1)/\psi(x)$. In $I_1$ and $I_3$, the change of variable $v = 1 + \frac{\psi(x)}{x}t$ yields

$$I_1 = \sqrt{\frac{\psi(x)}{x}} \int_{w_0}^\infty \frac{\bar{H}(x+t\psi(x))}{\bar{H}(x)} \frac{\mathrm{d}t}{(1+t(\psi(x)/x))\sqrt{1+(t/2)(\psi(x)/x)}\sqrt{2t}},$$

$$I_3 = 2\sqrt{\frac{\psi(x)}{x}} \int_0^\infty \frac{\bar{H}(x+t\psi(x))}{\bar{H}(x)} \frac{\mathrm{d}t}{(1+t(\psi(x)/x))\sqrt{1+(t/2)(\psi(x)/x)}\sqrt{2t}}.$$

In $I_2$, the change of variable $vy = x + t\psi(x)$ yields

$$I_2 = \frac{\psi(x)}{x} \int_{w_0}^\infty \frac{\bar{H}(x+t\psi(x))}{\bar{H}(x)} \frac{(y/x)\,\mathrm{d}t}{(1+t(\psi(x)/x))\sqrt{1-(y/x)^2+2t(\psi(x)/x)+(\psi^2(x)/x^2)t^2}}.$$

Let $J_i = \sqrt{x/\psi(x)}\,I_i$, $i = 1, 3$, and $J_2 = (x/\psi(x))I_2$. We start with $I_1$ and $I_3$. We will use the bound, valid for all $B, C > 0$,

$$0 \leq 1 - \frac{1}{(1+C)\sqrt{1+B}} \leq B/2 + C, \tag{22}$$

which follows from straightforward algebra and the concavity of the function $x \mapsto \sqrt{1+x}$. Applying this bound with $B = \frac{\psi(x)}{x}\frac{t}{2}$ and $C = \frac{\psi(x)}{x}t$ yields

$$0 \leq 1 - \frac{1}{(1+t(\psi(x)/x))\sqrt{1+(t/2)(\psi(x)/x)}} \leq \frac{5}{4}\frac{\psi(x)}{x}t.$$

We thus have

$$|J_1 - \sqrt{2\pi}\bar{\Phi}(\sqrt{2w_0})| + |J_3 - \sqrt{2\pi}|$$
$$\leq 3\chi(x)\int_0^\infty \Theta(t)\frac{\mathrm{d}t}{\sqrt{2t}} + \frac{15}{16}\sqrt{2\pi}\frac{\psi(x)}{x}.$$

Hence,

$$\frac{I_1}{I_3} = \bar{\Phi}(\sqrt{2w_0}) + O\left(\chi(x) + \frac{\psi(x)}{x}\right). \tag{23}$$

Now consider $J_2$. Applying the bound (22) with $C = t\psi(x)/x$ and

$$B = \left\{2t\frac{\psi(x)}{x} + \frac{\psi^2(x)}{x^2}t^2\right\} \Big/ \sqrt{1-(y/x)^2},$$

and making use of Hypothesis 1, we obtain

$$\left|J_2 - \frac{(y/x)\sqrt{2\pi}\varphi(\sqrt{2w_0})}{\sqrt{1-(y/x)^2}}\right|$$

$$\leq \frac{y/x}{\sqrt{1-(y/x)^2}}\chi(x)\int_0^\infty \Theta(t)\,\mathrm{d}t + \frac{y/x}{1-(y/x)^2}\frac{\psi(x)}{x}\left(\sqrt{1-\left(\frac{y}{x}\right)^2}+1+\frac{\psi(x)}{x}\right).$$

Choose $y = \rho x + \sqrt{1-\rho^2}\sqrt{x\psi(x)}z$ for some fixed $z \in \mathbb{R}$. Then, for large enough $x$, we have $0 < y < x$ and

$$\frac{y/x}{\sqrt{1-(y/x)^2}} = \frac{\rho}{\sqrt{1-\rho^2}} + O(\sqrt{\psi(x)/x}).$$

Thus,

$$\frac{I_2}{I_3} = \frac{\rho}{\sqrt{1-\rho^2}}\sqrt{\frac{\psi(x)}{x}}\varphi(\sqrt{2w_0}) + O\left(\chi(x) + \frac{\psi(x)}{x}\right). \tag{24}$$

For $z \geq 0$ and large enough $x$, plugging (23) and (24) into (20) yields

$$\theta(x,y) = \Phi(\sqrt{2w_0}) - \frac{\rho}{\sqrt{1-\rho^2}}\sqrt{\frac{\psi(x)}{x}}\varphi(\sqrt{2w_0}) + O\left(\chi(x) + \frac{\psi(x)}{x}\right).$$

For $z < 0$ and large enough $x$, plugging (23) and (24) into (21) yields

$$\theta(x,y) = \bar{\Phi}(\sqrt{2w_0}) - \frac{\rho}{\sqrt{1-\rho^2}}\sqrt{\frac{\psi(x)}{x}}\varphi(\sqrt{2w_0}) + O\left(\chi(x) + \frac{\psi(x)}{x}\right).$$

Now note that $w_0 = z^2/2 + O(\psi(x)/x)$, hence $\sqrt{2w_0} = |z| + O(\psi(x)/x)$. Thus, in both the cases $z \geq 0$ and $z < 0$, (4) holds. Let $z = z' + \rho\sqrt{\psi(x)/x}/\sqrt{1-\rho^2}$. A Taylor expansion of $\Phi$ and $\varphi$ around $z'$ yields (5). □

**Proof of Lemma 1.**

$$\frac{\bar{H}\{x+t\psi(x)\}}{\bar{H}(x)} - e^{-t} = \exp\left\{-\int_x^{x+t\psi(x)} \frac{ds}{\psi(s)}\right\} - e^{-t}$$

$$= \exp\left[-\int_0^t \frac{\psi(x)}{\psi\{x+s\psi(x)\}} ds\right] - e^{-t}.$$

Applying the inequality $|e^{-a} - e^{-b}| \leq |a-b|e^{-a \wedge b}$ valid for all $a, b \geq 0$ yields

$$\left|\frac{\bar{H}\{x+t\psi(x)\}}{\bar{H}(x)} - e^{-t}\right|$$

$$\leq \left|\int_0^t \frac{\psi\{x+s\psi(x)\} - \psi(x)}{\psi\{x+s\psi(x)\}} ds\right| \exp\left[-t \wedge \int_0^t \frac{\psi(x)}{\psi\{x+s\psi(x)\}} ds\right].$$

Let

$$\int_0^t \frac{|\psi\{x+s\psi(x)\} - \psi(x)|}{\psi\{x+s\psi(x)\}} ds = I(x,t)$$

and

$$\exp\left[-t \wedge \int_0^t \frac{\psi(x)}{\psi\{x+s\psi(x)\}} ds\right] = E(x,t).$$

*Case $\psi$ increasing.* If $\psi$ is non-decreasing, then $\psi' \geq 0$ and $\psi'$ is decreasing. Thus, for any $\delta > 0$ and large enough $x$, $\psi'(x) \leq \delta$ and

$$\int_0^t \frac{\psi(x)}{\psi\{x+s\psi(x)\}} ds \geq \int_0^t \frac{\psi(x)}{\psi(x) + s\psi'(x)\psi(x)} ds$$

$$= \int_0^t \frac{ds}{1+s\psi'(x)} \geq \int_0^t \frac{ds}{1+s\delta} = \frac{1}{\delta}\log(1+\delta t).$$

This implies that $E(x,t) \leq (1+\delta t)^{-1/\delta}$ for large enough $x$. Since $\psi$ is increasing and $\psi'$ is decreasing, we also have

$$I(x,t) \leq \psi'(x) \int_0^t \frac{s\psi(x)}{\psi\{x+s\psi(x)\}} ds \leq \psi'(x) \frac{t^2}{2}.$$

Thus, for any $\delta > 0$, $I(x,t)E(x,t) = O(|\psi'(x)|t^2(1+\delta t)^{-1/\delta})$.

*Case $\psi$ decreasing.* If $\psi$ is monotone non-increasing, then

$$\int_0^t \frac{\psi(x)}{\psi\{x+s\psi(x)\}} ds \geq t$$

and $E(x,t) \leq e^{-t}$. Also, since $|\psi'|$ is decreasing,

$$I(x,t) \leq |\psi'(x)| \int_0^t \frac{s\psi(x)}{\psi\{x+s\psi(x)\}} \, ds.$$

If $\psi$ has a positive limit at infinity, then $I(x,t) = O(\psi'(x)t^2)$. Otherwise, $\lim_{x\to\infty} \psi(x) = 0$ and (8) holds. This yields, for large enough $x$,

$$I(x,t) \leq |\psi'(x)| \int_0^t \frac{s\psi(x)}{\psi\{x+s\psi(x)\}} \, ds$$

$$\leq c_1|\psi'(x)| \int_0^t s e^{c_2 s\psi(x)} \, ds \leq c_1|\psi'(x)| \int_0^t s e^{s/2} \, ds.$$

Thus, $I(x,t) = O(|\psi'(x)|t^2 e^{t/2})$ and $I(x,t)E(x,t) = O(|\psi'(x)|t^2 e^{-t/2})$. This concludes the proof. $\square$

**Proof of Lemma 2.** Write

$$\left|\frac{\bar{H}\{x+t\psi(x)\}}{\bar{H}(x)} - e^{-t}\right| \leq \left|\frac{\bar{H}\{x+t\psi(x)\}}{\bar{H}\{x+t\psi_H(x)\}} - 1\right| e^{-t} \qquad (25)$$

$$+ \frac{\bar{H}\{x+t\psi(x)\}}{\bar{H}\{x+t\psi_H(x)\}} \left|\frac{\bar{H}\{x+t\psi_H(x)\}}{\bar{H}(x)} - e^{-t}\right|,$$

$$\frac{\bar{H}\{x+t\psi(x)\}}{\bar{H}\{x+t\psi_H(x)\}} = \exp\left[-\int_t^{t\psi(x)/\psi_H(x)} \frac{\psi_H(x)}{\psi\{x+s\psi_H(x)\}} \, ds\right]. \qquad (26)$$

If $\psi_H$ is increasing, then

$$\left|\int_t^{t\psi(x)/\psi_H(x)} \frac{\psi_H(x)}{\psi_H\{x+s\psi_H(x)\}} \, ds\right| \leq t\xi(x), \qquad (27)$$

$$\exp\left[-\int_t^{t\psi(x)/\psi_H(x)} \frac{\psi_H(x)}{\psi_H\{x+s\psi_H(x)\}} \, ds\right] \leq e^{t\xi(x)}. \qquad (28)$$

Since $\xi(x) \to 0$, gathering (27) and (28) yields, for large enough $x$,

$$\left|\frac{\bar{H}\{x+t\psi(x)\}}{\bar{H}\{x+t\psi_H(x)\}} - 1\right| e^{-t} \leq t\xi(x)e^{-t/2}. \qquad (29)$$

If $\psi_H$ is decreasing and $\lim_{x\to\infty} \psi_H(x) > 0$, then the ratio $\psi_H(x)/\psi_H\{x+s\psi(x)\}$ is bounded above and away from 0, so

$$\left|\int_t^{t\psi(x)/\psi_H(x)} \frac{\psi_H(x)}{\psi\{x+s\psi_H(x)\}} \, ds\right| \leq Ct\xi(x),$$

and (29) still holds.

If $\psi_H$ is decreasing and $\lim_{x\to\infty} \psi_H(x) = 0$, then applying (8) gives that the left-hand side of the previous equation is bounded by $Ct\xi(x)\exp\{2c_2\psi(x)t\}$. Thus, for large enough $x$, (29) still holds.

This provides a bound for the right-hand side of (25). The term (26) is bounded by Lemma 1. □

## Acknowledgements

We would like to thank two anonymous refereesfor useful remarks. We also acknowledge Free Software R. The first author gratefully acknowledges the support of the Natural Sciences and Engineering Research Council of Canada.